\documentstyle[12pt]{article}
\include{mssymb}
\newtheorem{prop}{Proposition}
\newtheorem{thm}{Theorem}
\newtheorem{cor}{Corollary}
\newtheorem{lemma}{Lemma}
\newtheorem{defn}{Definition}
\def\endproof{$\triangle$}
\def\P{{\Bbb P}}
\def\R{{\Bbb R}}
\def\C{{\Bbb C}}
\def\Z{{\Bbb Z}}
\def\A{{\Bbb A}}
\def\Zh{{\Bbb Z}[\frac{1}{2}]}
\def\Q{{\Bbb Q}}
\def\V{{\cal V}}
\def\W{{\cal W}}
\begin{document}
\title{Characteristic Classes and Quadric Bundles}
\author{Dan Edidin\thanks{Partially supported by an NSF postdoctoral
fellowship}  \mbox{ }and William Graham\thanks{Partially supported by an NSF
postdoctoral
fellowship}\\ Department of Mathematics\\ University of Chicago\\ Chicago IL
60637}
\date{}
\maketitle

\section{Introduction} \label{s.intro}
In this paper we construct Stiefel-Whitney and Euler classes in Chow
cohomology for algebraic vector bundles with nondegenerate quadratic
form.  These classes are not in the algebra generated by the Chern
classes of such bundles and are new characteristic classes in algebraic
geometry.  On complex varieties, they correspond to classes with the
same name pulled back from the cohomology of the classifying space
$BSO(N,\C)$.  The classes we construct are the only new characteristic
classes in algebraic geometry coming from the classical
groups (\cite{T2}, \cite{E-G-T}).

We begin by using the geometry of quadric bundles to study Chern
classes of maximal isotropic subbundles.  If $V\rightarrow X$ is
a vector bundle with quadratic form, and
if $E$ and $F$ are maximal isotropic subbundles of
$V$ then we prove (Theorem \ref{t.chern})
that $c_{i}(E)$ and $c_{i}(F)$ are equal mod 2.  Moreover, if the
rank of $V$ is $2n$, then $c_{n}(E) = \pm c_{n}(F)$, proving a
conjecture of Fulton.  We define Stiefel-Whitney and Euler
classes as Chow cohomology classes which pull back to Chern classes of
maximal isotropic subbundles of the pullback bundle.  Using the above
theorem we show (Theorem \ref{char.exist}) that these classes exist
and are unique, even
though $V$ need not have a maximal isotropic subbundle.  These
constructions also make it possible to give ``Schubert" presentations,
depending on a fixed maximal isotropic flag, of
the Chow rings of quadric and isotropic flag bundles.  The
presentation for the flag bundle is closely related to Fulton's work
on Schubert varieties in isotropic flag bundles (\cite{F1},
\cite{F2}).

An interesting aspect of this work is that it emphasizes
the difference between principal $SO(N)$ bundles in the Zariski and
\'etale topologies. As is well known (\cite{Sem-Chev}),
the group $SO(N)$, unlike $SL(N)$, is not
{\em special}; i.e. not all principal bundles which are
locally trivial in the \'etale topology are locally
trivial in the Zariski topology. This fact manifests itself
as follows. The Stiefel-Whitney classes in $A^*(X;\Z/2\Z)$
we construct are for bundles which are locally trivial
in the Zariski topology. In the
\'etale topology it is not always possible to construct these
classes, as there are bundles for which the corresponding
topological Stiefel-Whitney classes are not algebraic
\cite{T1}. Similarly, we construct an integral Euler
class for Zariski locally trivial $SO(2n)$ bundles, but for
general $SO(2n)$ bundles only a power of 2 times the Euler
class is integral \cite{T1}.

The paper is organized as follows. Section
\ref{s.intro} is the introduction. Section \ref{s.prelim}
is largely a collection of standard
and/or easy facts
about vector bundles with quadratic form, quadric bundles
and isotropic flag bundles which we use and which
have not been adequately presented in the literature.

The remainder of the paper is new. In Section \ref{s.chern} we prove
the theorem on Chern classes of maximal isotropic subbundles,
including Fulton's conjecture.  Section \ref{s:Euler} contains the
construction of Stiefel-Whitney and Euler classes for bundles with
quadratic form which are locally trivial in the Zariski topology. In
Section \ref{s:topology} we show that for smooth complex varieties the
classes we construct map to the classes with the same names in
singular cohomology.

Section \ref{s.halfchow}
deals with bundles whose structure group reduces to $SO(2n)$ but which
are not locally trivial. In particular, we show that there is a
characteristic class $y_n \in A^nX$ such that $\frac{1}{2^{n-1}}y_n$
is an Euler class.  We also give a natural (with $\Zh$ coefficients)
presentation of the Chow ring of the quadric bundle as an algebra with
two generators and two relations, such that the Euler and Chern classes
of $V$ are the coefficients in the relations.
This presentation
is the $SO(2n)$ analogue of the presentation of $A^*\P(V)$ over
$A^*X$,
where the Chern classes of V are the coefficients in the sole relation.
Viewing the isotropic flag variety $Fl(V)$ as a tower of quadrics bundles,
we also describe $A^*(Fl(V);\Zh)$. Finally, in Section \ref{s.zchow}
we give a ``Schubert'' presenations for $A^*Q$ and $A^*(Fl(V))$, which
depend on a fixed maximal isotropic subbundle $F \subset V$.

{\bf Acknowledgements.} We would like to thank William Fulton for his
generous help, in particular for informing us of his conjecture and
suggesting that the Euler class was related to top Chern classes of
maximal isotropic subbundles. Some of the facts presented in
Section 2 are based on his lectures at the University of Chicago.
We are also grateful to Burt Totaro for his example showing that
the Euler class need not be integral for bundles which are
not locally trivial in the Zariski topology.

\section{Preliminaries} \label{s.prelim}
Throughout this paper, unless otherwise noted, all schemes are assumed to be
of finite type over an arbitrary field.

\medskip

{\bf Operational Chow groups} If $X$ is a scheme, then $A_*X$ denotes
the Chow homology groups defined in \cite[Chapter 1]{Fulton}, and
$A^*X$ is the operational Chow cohomology ring of \cite[Chapter
17]{Fulton}.  An element $c \in A^pX$ is a collection of homomorphisms
$A_*X' \rightarrow A_{*-p} X'$ for all $X' \rightarrow X$ which is
compatible with proper pushforward, flat pull-back, and interstions.
By construction, if $X \stackrel{f} \rightarrow Y$ is any map there is
always a ring homomorphism $f^*:A^*Y \rightarrow A^*X$.  If $X$ is
regular, then $A^*X$ is just the usual Chow groups with multiplication
given by intersection product.  Now if $R$ is any coefficient ring, we
can also define $A^*(X;R)$ as maps from $A_*X' \otimes R
\rightarrow A_*X' \otimes R$ for all $X' \rightarrow X$
satisfying the same conditions as above. There is always
a homomorphism $A^*X \otimes R \rightarrow A^*(X;R)$. If
$X$ admits a resolution of singularities then standard arguments
(cf. \cite{Kimura}) show that $X$ is an isomorphism.

\medskip

{\bf Quadratic Forms} A quadratic form $q$ on a vector space $V$ will
always mean a {\em hyperbolic} quadratic form: in other words, we
assume that $V$ (of dimension $2n$ or $2n+1$) has an isotropic
subspace of dimension $n$ (if the ground field is algebraically
closed of characteristic not equal to 2,
this is automatically satisfied).  We can choose a basis
$v_{1}, v_{2}, \ldots$ of $V$ such that the quadratic norm $q(v,v)$ of
a vector $v = \sum x_{i}v_{i}$ is given by $$q(v,v) = x_{1}x_{2n} +
\ldots + x_{n}x_{n+1} \mbox{ (if dim $V = 2n$)}$$ $$q(v,v) =
x_{1}x_{2n+1} + \ldots + x_{n}x_{n+2} + a x_{n+1}^2 \mbox{ (if dim $V
= 2n+1$)}$$ where $a$ is a nonzero scalar.  The reason our form must
be hyperbolic is that otherwise there may be no nonzero isotropic
vectors in $V$.  For example, there are no nonzero isotropic vectors
for $\R^2$ with quadratic form $x_{1}^{2} + x_{2}^{2}$.  While the
theorems we prove are for schemes of arbitrary fields, we will abuse
notation and refer to the group that preserves a quadratic form on a
$N$ dimensional vector space as $O(N)$, even though this notation is
only correct when the field is algebraically closed.

Let $\mbox{Iso}(V)$ be the variety of maximal isotropic subspaces
of $V$.  The groups $SO(N)$ and $O(N)$ act on $\mbox{Iso}(V)$.  The
$O(N)$ action is always transitive, but the $SO(N)$ action is only
transitive if $N$ is odd. If $N$ is even, there are two $SO(N)$
orbits and two maximal isotropic subspaces $E$ and $F$ of $V$ lie in the
same $SO(N)$ orbit if and only if $\mbox{dim }E \cap F
\equiv n(\mbox{mod }2)$.  We will say that $E$ and $F$ are in the same
family if they are in the same $SO(N)$ orbit, and in opposite families
otherwise.

Let $Fl(V)$ denote the flag variety of complete isotropic flags of a
certain length in $V$.  This length is $n$ if dim $V=2n+1$ and $n-1$
if dim $V=2n$.  If dim $V=2n$ we let $Fl_n(V)$ denote the variety of
isotropic flags of length $n$.  A point in one of these varieties
corresponds to a flag $E_1 \subset E_2 \subset \ldots $ of isotropic
subspaces of $V$ (where $\mbox{dim }E_i=i$).  In the even dimensional
case, every length $n-1$ isotropic flag can be completed to exactly
two isotropic length $n$ flags, and the space $Fl_n(V)$ is a
disconnected double cover of $V$.  Again, the groups $SO(N)$ and $O(N)$
act. The $SO(N)$ action on $Fl(V)$ is transitive, but there are two
$SO(2n)$ orbits on $Fl_n(V)$; two length $n$ isotropic flags $E.$ and
$F.$ are in the same $SO(N)$ orbit if and only if $\mbox{dim }E_{n}
\cap F_{n}
\equiv n(\mbox{mod }2)$.

Let $Q \subset \P(V)$\footnote{Throught this paper
$\P(V)$ will denote the scheme of lines in $V$.}
denote the quadric of isotropic vectors.  If $E$
and $F$ are maximal isotropic subspaces then $\P(E)$ and $\P(F)$ are
subvarieties of $Q$ called rulings.  $\P(E)$ and $\P(F)$ represent the same
class in $A_{n-1}(Q)$
if and only if $E$ and $F$ lie in the same $SO(N)$ orbit.  Thus, if
dim $V=2n$, there are two families of rulings giving
rise to two Chow classes $e$ and $f$, and, letting $h$ denote the
hyperplane section, we have
$h^{n-1} = e+f$.  If dim $V=2n+1$ there is only one family
and $h^{n}=2e$.

\begin{lemma} \label{quad}
$(a)$ If $N=2n$ is even, then a basis for $A_*(Q)$ is $1,h,h^2, \ldots,
h^{n-1},e,he,\ldots,h^{n-1}e$.

$(b)$ If $N=2n+1$ is odd, then a basis for $A_*(Q)$ is $1,h,h^2, \ldots,
h^{n-1},e,he,\ldots,h^{n-1}e$. \endproof
\end{lemma}

Remark: Because a quadric has a decomposition into affine cells,
$A_*(U \times Q) = A_*(U) \otimes A_*(Q)$. This fact will be essential
when we compute the Chow rings of quadric bundles.

\medskip

Proof: (a) In any characteristic the hyperbolic quadric has an affine
cellular decompostion with 1 cell in every dimension, except 2 cells
in dimension $n-1$. Hence $A_i(Q) = \Z$ for $i \neq n-1$, and
$A_{n-1}(Q)=\Z \oplus \Z$.  Since $e$ is a class of a linear space,
$h^{n-1}e= 1$. Hence, for $i \leq n-1$, $h^i$ is a generator of
$A_{n-i}(Q)$. Likewise, the $h^ie$ are also generators of $A_{i}(Q)$
for $i \leq n-1$.  This proves (a).  The proof of (b) is similar.
\endproof

\medskip

In the sequel we use the following relations in the Chow ring of
an even-dimensional quadric.
\begin{lemma} \label{rels}
Let $Q \subset \P^{2n-1}$ be an even-dimensional quadric.
Then

$e^2 = f^2=0$ and $ef = 1$ if $n$ is even, and

$e^2 = f^2=1$ and $ef = 0$ if $n$ is odd.
\end{lemma}
Proof: First note that $e$ and $f$ both have degree 1,
so $h^{n-1}e=h^{n-1}f = 1$. Since $h^{n-1} = e+f$, it follows that
$e^2 = f^2$. Also, since $h^{2n-2} = (e+f)^2=2$, it follows either
that $e^2=f^2=1$ and $ef=0$, or, $e^2=f^2=0$ and $ef=1$.  The linear
system spanned by $e$ is an $SO(2n)$ orbit of
isotropic $n$-planes.  One such $n-$plane is given by $x_1=x_2 =
\ldots =x_n=0$. When $n$ is even, this plane can be taken to the plane
$x_{n+1} = x_{n+2} = \ldots = x_{2n}=0$ via an element of $SO(2n)$ (an
even permutation).  Projectivizing, we see that there are two disjoint
elements in the linear system spanned by $e$. Hence $e^2=0$.

When $n$ is odd, the same argument shows that
there are elements in the linear system corresponding to
$e$ that are disjoint from elements of $f$. Thus $ef=0$.
\endproof

\medskip

The constructions above also work for vector bundles.  We adopt the
convention that a quadratic form $q$ on a vector bundle $V \rightarrow
X$ means a section of $S^{2}(V^*)$ which restricts to a nondegenerate
hyperbolic quadratic form on each fiber.  The constructions of $Q$,
$Fl(V)$, and $Fl_n(V)$ all carry over to the bundle case.  A maximal
isotropic subbundle is defined to be an isotropic subbundle of rank
$n$ (where the rank of $V$ is $2n$ or $2n+1$).  If the rank of $V$ is
even and the structure group reduces to $SO(2n)$ (this is
automatically satisfied in the Zariski locally trivial case -- see
below) then the pullback of $V$ to $Fl(V)$ has two tautological
maximal isotropic subbundles.  If the rank is odd then there is only
one tautological maximal isotropic subbundle over $Fl(V)$.  These
facts will be used extensively in the sequel.

The description of $SO(2n)$
orbits of maximal isotropic subspaces yields the following fact about
maximal isotropic subbundles.

\begin{prop} \label{p.const}
If $E$ and $F$ are maximal isotropic subbundles of $V \rightarrow
X$, where the rank of $V$ is even, then $\mbox{{\em dim }} E_x \cap
F_x$ is constant {\em mod} $2$ on each component of $X$.
\endproof \end{prop}

We say the pair
$(V,q)$ is Zariski locally trivial (or simply locally trivial)
if there exists a Zariski open covering $\{U_i\}$ of $X$, such that
$V_{|U_i} \simeq X \times \A^{N}$ and the under the isomorphisms
the quadratic norm on a fiber over $u \in U_i$
is given by:
$$\mbox{(rank  $V=2n$) } q(v,v) = v_1v_{2n} + \ldots + v_nv_{n+1}$$
$$\mbox{(rank  $V=2n+1$) } q(v,v) = v_1v_{2n+1} + \ldots + v_{n}v_{n+2} +
f(u)v_{n+1}^2.$$
Here $v=(v_1, \ldots , v_N)$ is a vector in the fiber over $u$, and
$f$ is a nowhere vanishing function on $U_i$.
(Swan \cite[Cor 1.2]{Swan} has shown that
any $(V,q)$ is locally trivial in the \'etale topology; i.e.,
there is an \'etale cover $\{U_i \stackrel{\pi_i}
\rightarrow X \}$ such that $(\pi_i^*V,\pi_i^*q)$ is as above).

\medskip

{\bf Remarks:}
Not all bundles with quadratic form are locally trivial in
the Zariski topology, as the following example shows.
{\bf Example:}
Let $X=\C^*$ be the $1$ dimensional torus. Let $V = X \times \C^2$
be a trivial rank two vector bundle. Define a quadratic form
on $V$ by the rule that over a point $\lambda\in \C^*$
the quadratic norm is $v_1^2 + \lambda v_2^2$.  The pair $(V,q)$ is
locally trivial in the \'etale topology but not in the Zariski topology.

Maximal isotropic subbundles need not exist, even if $V$
is Zariski locally trivial; we give a simple example in Section
\ref{s.halfchow}.  However, if $V$ has a maximal
isotropic subbundle, then the pair $(V,q)$ is locally trivial.  If the
rank of $V$ is odd, then it is not always possible to
choose trivializations so that the functions
$f(u)$ above are identically $1$,
even if $V$ has a maximal isotropic subbundle.
For instance, if $V$ is a self-dual
line bundle, then the trivializations can be so chosen if and only if
$V$ is the trivial line bundle.

If the pair $(V,q)$ is locally trivial in the Zariski
topology, then so is the associated quadric bundle.  This assertion is
obvious if the rank is even.  In the odd rank case, it suffices to
show that if $V \rightarrow U$ is a trivial vector bundle with
quadratic form as above, then the quadric bundle
$Q$ on $U$ is isomorphic to $U \times Q_{pt}$, where $Q_{pt}$ is the
quadric over a point.  The
desired isomorphism maps
$u \times (v_{1},\ldots,v_{2n+1})$ to
$(f(u)v_{1},\ldots,f(u)v_{n},v_{n+1},\ldots, v_{2n+1})$.

Note that if $(V,q)$ is Zariski locally trivial of rank $2n$
then the structure group of $V$ reduces to $SO(2n)$.
The reason is that the principal bundle is also Zariski locally
trivial.  This bundle is a priori an $O(2n)$ bundle,
but any Zariski locally trivial bundle on a connected space
has the same number of components as the fiber.  Thus (assuming $X$
connected) the $O(2n)$ principal bundle has two components, implying
that the structure group reduces.
If the rank is $2n+1$ we cannot conclude that the structure group reduces to
$SO(2n+1)$, because the associated $O(2n+1)$ principal bundle need not
be Zariski locally
trivial (since the function $f(u)$ need not be identically 1).

The following lemma well be needed in the sequel.
\begin{lemma} \label{thom.bundle}
Let $f:Y \rightarrow X$ be a bundle with complete fibers
which is locally trivial in the Zariski topology.
Then
$f^*$ is injective for Chow cohomology with any coefficients.
\end{lemma}
Proof: Because $Y$ is a locally
trivial in Zariski topology, $f:Y \rightarrow X$
is a proper Chow envelope in the sense of \cite[Chapter 18]{Fulton}.
The lemma follows from \cite[Lemma 2.1]{Kimura}. \endproof

\medskip

{\bf Remark on characteristic 2.} The results of this paper
hold in characteristic 2 (and if we assume that the constant $a$
in the odd rank quadratic form is a square in the field of definition).
The reason is that the
integral quadrics $x_1x_{2n} + \ldots x_{n-1}x_{n+1}$ or
$x_1x_{2n+1} + \ldots x_{n}x_{n+2} + x_{n+1}^2$ are smooth over
$Spec\;\Z$. Likewise we can construct an isotropic flag variety which
is smooth over $Spec\;\Z$. By specializing from
characteristic 0, we can prove all the necessary intersection
theoretic properties of the quadric and flag varieties defined
over ${\Bbb F}_2$. Since our constructions depend only on intersection
theory, the results follow. For consistency of notation
we take $SO(N)$ in characteristic 2 to mean the connected
component of the identity in $O(N)$.

Note that in characteristic 2, even if the field is algebraically
closed not all forms are hyperbolic. For example the quadratic form
$x_1^2 + x_2^2$ over $\overline{{\Bbb F}}_2$ is not hyperbolic.

\section{Chern Classes of Isotropic Subbundles} \label{s.chern}
The purpose of this section is to prove the following theorem.
Part (c) is Fulton's original conjecture.

\begin{thm}  \label{t.chern}
Let $V$ be a vector bundle on a connected scheme $X$, of
rank $N$ equal to $2n$ or $2n+1$, equipped with
a non-degenerate quadratic form, and let $d=N-n$.  Let $E$
and $F$ be maximal isotropic subbundles of $V$.

$(a)$ There exist
classes $c_{i}$ and $d_{i}$ in $A^*(X)$ such that $2c_{i} =
c_i(E)+ c_i(F)$ and $2d_{i} = c_i(V/E) +c_i(V/F)$.  These classes are
functorial for maps $X^\prime \rightarrow X$.

$(b)$ If $E$ and $F$ are in the same (resp. opposite) families, then
$c_d(V/E) = c_d(V/F)$ (resp. $-c_d(V/F)$).

$(c)$ If the rank of $V$ is $2n$, then $c_n(E)=c_n(F)$ if
$E$ and $F$ are in the same family, and $c_n(E)=-c_n(F)$ if $E$ and
$F$ are in opposite families.
\end{thm}
Proof: Note that (c) follows immediately from (b), because in the even
rank case, we may identify $V/E$ and $V/F$ with $E^*$ and $F^*$ respectively.

If $Y \rightarrow X$, then $c_i(V/E) \cap [Y]$ and
$c_i(V/F) \cap [Y]$ are both supported in $A_*(Y)$. Replacing $X$ with
$Y$ it suffices to check the relations after capping with $[X]$.
In the remainder of the proof, we will use the notations $c_i(V/E)$, etc.,to
mean $c_i(V/E) \cap [X] \in A_*(X)$.
Let $Q\rightarrow X$ be the quadric bundle associated to the quadratic
form on $V$. Let $i:Q \hookrightarrow \P(V)$ be the inclusion, and
let $\pi:\P(V) \rightarrow X$ and $\rho:Q \rightarrow X$ be the projections.
Because $E$ and $F$ are isotropic subbundles,
$\P(E)$ and $\P(F)$ are subvarieties of $Q$.

Now we will examine $[\P(E)]$ and $[\P(F)]$ in the Chow groups of
$Q$. Set $\gamma = [\P(E)]$, and let $h$ be the pullback to $A_*(Q)$
of the hyperplane class $H$ of $\P(V)\rightarrow X$.
Since $\P(E) \subset Q$ is a regular embedding
we can define a bivariant class $\gamma \in A^*(Q \rightarrow X)$
by the formula
$\alpha \rightarrow \gamma \cap \rho^*\alpha \in A_*(Q)$.

\begin{lemma} \label{gen}
Every $\beta \in A_*(Q)$ can be written as
$$\sum_{i=0}^{i=n-1} h^i \rho^*(\alpha_i) + \sum_{i=0}^{i=n-1}
h^i \gamma \cap \rho^*(\beta_i).$$
\end{lemma}

Proof of Lemma \ref{gen}: By Lemma \ref{quad}, the classes listed
generate the Chow groups of the fibers. Because $V$ has a maximal
isotropic subbundle, $Q \rightarrow X$ is locally trivial in
the Zariski topology.
The conclusion then follows by
Noetherian induction as in \cite[Proposition 1.9]{Fulton}. \endproof

\medskip

We now finish the proof of the theorem.  First suppose that the
$E$ and $F$ are in the same family (this is automatically satisfied in
the odd rank case).  Then $[\P(E)]$ and $[\P(F)]$ agree on each fiber
of $Q$. Thus in $A_*(Q)$, $$[\P(E)] -[\P(F)]= \sum_{i > 0} h^{d-1-i}
\cap \rho^*(\alpha_{i}),$$ for some classes $\alpha_{i}$ in $A_*(X)$.
But
$$i_*(h^j \cap \rho^* \alpha) = c_1({\cal O}(Q)) \cap H^j \cap
\pi^*(\alpha).$$
Since the quadric bundle corresponds to a section of
${\cal O}_{\P(V)}(2)$, $c_1({\cal O}(Q))=2H$.  Thus
$$i_*([\P(E)] -[\P(F)])= \sum_{i>0} 2H^i \cap \pi^*(\alpha_{d-i}).$$
On the other
hand, \cite[Example 3.2.17]{Fulton} implies that
$$i_*([\P(E)] - [\P(F)])=
\sum_{i \geq 0} H^i \cap \pi^* (c_{d-i}(V/E) - c_{d-i}(V/F)).$$Comparing
these expressions, we see that $c_d(V/E) = c_d(V/F)$, which proves
part (b).  Moreover, defining $d_{i} = \alpha_{i} + c_i(V/F)$, we see that
$2d_{i} = c_i(V/E) +c_i(V/F)$ as desired.
In the even case, to get the classes $c_i$, we simply
define $c_i = (-1)^i d_{i}$.  In the odd case, we no longer have $V/E$
and $V/F$ isomorphic to $E^*$ and $F^*$.  However, an easy computation
(using the fact that $c_1(V) = c_1(E^{\perp}/E) = c_1(F^{\perp}/F)$)
shows that
$$
c(E^*) + c(F^*) = \frac{c(V/E) + c(V/F)}{1+c_1(V)}.
$$
Therefore, if we define the $c_i$ by the formula
 $$
\sum_{i=0}^{n} (-1)^i c_i = \frac{\sum_{i=0}^{n} d_i}{1+c_1(V)}.
$$
we get $2c_{i} = c_i(E) +c_i(F)$, as desired.  The classes $c_i$ and
$d_i$ are obviously functorial for maps $X^\prime \rightarrow X$.

This proves the theorem in the case where $E$ and $F$ are in the same
family.  In the case where $E$ and $F$ are in opposite families we
compare $[\P(E)]$ and $h^{n-1} - [\P(F)]$ to obtain the result.
\endproof

\section{Construction of Euler and Stiefel-Whitney Classes}
\label{s:Euler}

Let $V \rightarrow X$ be a vector bundle of rank $N=2n$ or $N=2n+1$
with a nondegenerate
quadratic form $q:V \otimes V \rightarrow {\cal O}$.

\begin{defn}
(a) A {\em Stiefel-Whitney class} for $V \rightarrow X$ is
a class $w_{2i} \in A^i(X;\Z/2\Z)$
such that for any $f:Y \rightarrow X$
and any maximal isotropic subbundle $E \subset f^*V$ we have
$f^*w_{2i} \equiv c_i(E)(\mbox{mod } 2)$

(b) If the rank of $V$ equals $2n$, then an {\em Euler class}
for $V \rightarrow X$ is a class $x \in A^n(X)$ such that
for any $f:Y \rightarrow X$, and any maximal isotropic
subbundle $E \subset f^*V$ we have $f^*x = \pm c_n(E)$.
\end{defn}
The purpose of this section is to prove the existence
of these classes when the pair $(V,q)$ is locally trivial in
the Zariski topology.
For the remainder of this section we will therefore assume
that $(V,q)$ is Zariski locally trivial.
We may repeat this assumption for emphasis.

{\bf Remark.} Note that if $E \subset V$ is a maximal isotropic
subbundle, then $c_n(E)$ is an Euler class for $V$.
In Section \ref{s:topology} we will explain the connnection
between our defintions and those given in topology.

The next proposition states some basic properties
of Stiefel-Whitney and Euler classes.

\begin{prop} \label{p.basic}

$(1)$ Pullbacks of Stiefel-Whitney (resp. Euler) classes are
Stiefel-Whitney (resp. Euler) classes for the pullback bundle.

$(2)$ Stiefel-Whitney classes are unique.

$(3)$ An Euler class is unique up to sign; i.e., if
$x,y$ are both Euler classes for a bundle $V_n$, then
$x = \pm y$.

$(4)$ If $x$ is an Euler class of $V_n$, then
$x^2 = (-1)^{n}c_{2n}(V_n)$.
\end{prop}
Proof: Property (1) follows from the definition.  Next
consider the map $f:Fl(V) \rightarrow X$.  The bundle
$f^*V$ has a maximal isotropic subbundle $E$.
Suppose $w_{2i}$ and $w_{2i}^\prime$ are Stiefel-Whitney
classes in degree $A^i(X)$.
By definition,
$$f^*w_{2i}=f^*w_{2i}^\prime =c_i(E)(\mbox{mod 2}).$$
Since $Fl(V) \rightarrow X$ is locally trivial, $f^*$
is injective (Lemma \ref{thom.bundle}), so
$w_{2i} = w_{2i}^\prime$ in $A^*(X;\Z/2\Z)$, proving
$(2)$.
Likewise,
$f^*x$ and $f^*y$ are equal to $\pm
c_n(V^\prime_n)$, so
$x=\pm y$ as Chow cohomology classes. This proves (3).

Finally, note that, via the quadratic form, we can identify
$V_n^{\prime *} = f^*V_n/V_n^\prime$. Thus,
$$c_{2n}(f^*V_n) = c_n(V_n^\prime) c_n(V_n^{\prime *})
= (-1)^nf^*(x^2).$$
Since $f^*$ is injective, (4) follows. \endproof

\begin{thm} \label{char.exist}
Let $V \rightarrow X$ be a vector bundle
with a nondegenerate quadratic form $q$.
If the pair $(V,q)$ is locally trivial in the Zariski topology
then $V$ has Stiefel-Whitney classes $w_{2i} \in A^i(X;\Z/2\Z)$.
If the rank of $V$ is $2n$, then $V$ has two Euler classes
$x_n$ and $-x_n$ in $A^n(X)$.
\end{thm}
{\bf Remark 1.} Proposition \ref{p.basic} and Theorem \ref{char.exist}
together show that Stiefel-Whitney and Euler classes are characterstic
classes for vector bundles with quadratic form which are locally trivial
in the Zariski topology.

\medskip

{\bf Remark 2.} In topology the sign of the Euler class is determined
by the choice of an orientation of the bundle, or in other words, a
reduction of structure group to $SO(N)$.  The same holds in the
algebraic case.  This is equivalent to choosing a maximal isotropic
subbundle of the pullback of $V$ to the flag bundle.
\medskip

Proof: Consider the Cartesian diagram
$$\begin{array}{ccc}
Fl(V) \times_{X} Fl(V) & \stackrel{g'}\rightarrow & Fl(V) \\
\downarrow\scriptsize{f'} & & \scriptsize{f}\downarrow\\
Fl(V) & \stackrel{g}\rightarrow & X
\end{array}$$
Let $E \subset f^*V$ be one of the two tautological maximal
isotropic subbundles. (The choice of $E$ is equivalent to choosing an
orientation.)  To show that $w_{2i}$ exists we must
show that $c_i(E)(\mbox{mod 2})$ is the pullback (mod 2)
of a class $A^*X$.

Let $F \subset g^*V$ be a maximal isotropic subbundle. Then $c_i(F)
(\mbox{mod 2})$ is a Stiefel-Whitney class for $g^*V$.
By Theorem \ref{t.chern},
$$f^{\prime*}c_i(F) \equiv c_{i}(g^{\prime *}E) \equiv
g^{\prime *} c_{i}(E)(\mbox{mod 2}).$$
Now $Fl(V) \rightarrow X$ is a proper Chow envelope, so
by \cite[Theorem 2.3]{Kimura} $c_i(E)$ is a pullback
(mod 2).

To complete the proof we must show that for any
$q:X^\prime \rightarrow X$, and any maximal isotropic
subbundle $F \subset q^*V$, we have $q^*w_{2i} \equiv c_i(F)
(\mbox{mod 2})$.
To prove this consider the Cartesian diagram
$$\begin{array}{ccc}
Fl(q^*V) & \stackrel{q'}\rightarrow & Fl(V) \\
\downarrow\scriptsize{f'} & & \scriptsize{f}\downarrow\\
X' & \stackrel{g}\rightarrow & X
\end{array}$$
Chasing the diagram, and applying Theorem \ref{t.chern}, it
follows that $f^{\prime *}q^*(x_n) = \pm f^{\prime *}(c_n(F))$. Since
$f^*$ is injective mod 2, the existence of Stiefel-Whitney
classes follows.

If $V$ has rank $2n$, then we can prove the existence of
Euler classes in $A^n(X)$ by the method above, using the
fact that $c_n(E)=\pm c_n(F)$ for any two maximal
isotropic subbundles.
\endproof

{\bf Remark.} If the bundle and quadratic form are not locally trivial
in the Zariski topology  it is still possible to construct (using our
techniques)
classes $\pm y_n \in A^n(X)$
(Theorem \ref{euler.haexist}) such
that $\pm \frac{y_n}{2^{n-1}}$ are Euler classes.
On the other hand, the calculations of B. Totaro
\cite{T1} show that Stiefel-Whitney classes do not exist
and that the Euler class is not integral for all quadratic vector bundles
which are not locally trivial.
In particular, he shows that on a variety
approximating the classifying space $BSO(4)$ there is a
tautological quadratic vector bundle where the topological
Stiefel-Whitney and Euler classes are not
represented by complex manifolds.

\section{Connections with topology}
\label{s:topology}

In this section we explain the connection between the definitions of
Euler and Stiefel-Whitney classes in the algebraic case and the usual
definitions of these classes in topology.  One consequence of this
connection is a proof of a topological analogue of Theorem \ref{t.chern}.
Throughout this
section, unless otherwise stated, we work in the setting of topology
rather than algebraic geometry.

We begin by recalling a few facts about characteristic classes and
classifying spaces.  The space $BO(N,\C)$ is the classifying space for
rank $N$ complex vector bundles with nondegenerate quadratic form, and
$BO(N,\R)$ is the classifying space for rank $N$ real vector bundles.
Let $\V \rightarrow BO(N,\C)$ and $\W \rightarrow BO(N,\R)$ be the
universal vector bundles.  The bundle $\V$ (resp. $\W$) has a
nondegenerate (resp. positive definite real) quadratic form.
Complexifying the quadratic form on $\W$ gives a nondegenerate
quadratic form on $\W \otimes \C$.  Because topologically the groups
$O(N,\C)$ and $O(N,\R)$ are homotopy equivalent, we can identify the
spaces $BO(N,\C)$ and $BO(N,\R)$.  Under this identification, $\V = \W
\otimes \C$, as bundles with quadratic form.  Because this is the
universal case, we have the following lemma.

\begin{lemma}
\label{l:class}

Let $V \rightarrow X$ be a complex vector bundle with nondegenerate
quadratic form.  Then there exists a real vector bundle $W \rightarrow
X$ with positive definite quadratic form, such that $V \simeq W
\otimes \C$ as vector bundles with quadratic form.  Moreover, the
classifying maps $X \rightarrow BO(N,\C)$ of $V$ and $W$ coincide.
\endproof
\end{lemma}

The cohomology ring $H^{*}(BO(N,\C); \Z / 2 \Z)$ is isomorphic to the
polynomial ring \\ $\Z / 2 \Z [w_{1}, \ldots, w_{N}]$.
The Stiefel-Whitney classes of a complex vector
bundle $V$ with nondegenerate quadratic form (resp. real vector bundle
$W$) over $X$ are defined to be the pullbacks of these classes via the
classifying map to $H^{*}(X; \Z / 2 \Z)$.  If the structure group of
$V$ (resp. $W$) reduces to $SO(N,\C)$ (resp. $SO(N,\R)$), then the
classifying map $X \rightarrow BO(N,\C)$ lifts (in exactly two ways,
assuming $X$ connected) to a map $X \rightarrow BSO(N,\C)$.  If
$N=2n$, we single out a universal class $x_{n} \in H^{n}(BSO(2n,\C);
\Z)$; the Euler classes of $V$ (resp. $W$) are defined to be the
pullbacks of $x_{n}$ to $H^{n}(X; \Z)$ via the two lifts.  These
pullback classes differ only by sign.

If $V$ and $W$ are as in Lemma \ref{l:class}, then because the
classifying maps of $V$ and $W$ coincide, their Stiefel-Whitney and
Euler classes coincide as well.

\begin{prop}
\label{p:topiso}
Let $V \rightarrow X$ be a complex vector bundle with nondegenerate
quadratic form and let $E$ be a maximal isotropic subbundle of $V$.
Then the even Stiefel-Whitney classes of $V$ are the mod $2$
reductions of the Chern classes of $E$.  Moreover, if the rank of $V$
is $2n$, then the structure group of $V$ reduces to $SO(2n,\C)$, the
Euler classes of $V$ are $\pm c_{n}(E)$, and the odd Stiefel-Whitney
classes of $V$ vanish.

\end{prop}

Proof: Choose $W$ as in Lemma \ref{l:class}.  Let $V_{\R}$ and
$E_{\R}$ denote $V$ and $E$ viewed as real vector bundles.  We have
maps of real bundles
$$
E_{\R} \hookrightarrow V_{\R} \simeq W \oplus iW \rightarrow W .
$$
The composition $E_{\R} \rightarrow W$ is
injective, since $E$ is isotropic and the quadratic form on $W$ is
positive definite.  Assume now that the rank of $V$ is $2n$.  Then the
ranks of $E_{\R}$ and $W$ are equal, so $E_{\R} \simeq W$.  Since $W$
is isomorphic to the realification of a complex vector bundle, it is
orientable as a real vector bundle, which means that the structure
group of $W$ reduces to $SO(2n,\R)$.  This implies that the structure
group of $V$ reduces to $SO(2n,\C)$.  The Stiefel-Whitney and
Euler classes of $V$ equal those of $W$, which in turn equal those of
$E_{\R}$.  The assertions of the proposition then follow from
\cite[Problem 14-B and Definition, p. 158]{Milnor-Stasheff}.

This proves the proposition in the even rank case.  In the odd rank
case, we can write $W \simeq E_{\R} \oplus L$, where $L$ is a real
line bundle.  The proposition then follows from the Whitney product
formula and \cite[Problem 14-B]{Milnor-Stasheff}.
\endproof

\medskip

The preceding proposition shows that the Euler and Stiefel-Whitney
classes we have defined in Chow cohomology bear the same relation to
the Chern classes in Chow cohomology as the Euler and Stiefel-Whitney
classes in topology do to the Chern classes in topology, and hence
justifies our use of these names.

In the topological setting, the properties of Stiefel-Whitney and
Euler classes given in Proposition \ref{p:topiso} yield the
following corollary, which is a topological analogue of
Theorem \ref{t.chern}.  (In the algebraic case, the logic is reversed:
we use Theorem \ref{t.chern} to prove the existence of Stiefel-Whitney and
Euler classes.)

\begin{cor}
\label{c:weakeuler}

Let $V \rightarrow X$ be a complex vector bundle with
nondegenerate quadratic form, of rank $N$ equal to $2n$ or $2n+1$,
and suppose that $E$ and $F$ are maximal isotropic subbundles of $V$.
Then the mod $2$ reductions of $c(E)$ and $c(F)$ are equal.  Moreover,
if $N=2n$, then $c_{n}(E) = \pm c_{n}(F)$.
\endproof
\end{cor}

For smooth complex varieties the existence of a natural map from
Chow cohomology to singular cohomology yields the following
proposition.
\begin{prop}
\label{p:correspond}
Let $X$ be a smooth complex variety and let $V \rightarrow X$ be a vector
bundle with quadratic form, locally trivial in the Zariski topology.
Then, under the natural map $co:A^i(X;R) \rightarrow H^{2i}(X;R)$
(for $R=\Z/2\Z$ or $R=\Z$)
the algebraic Stiefel-Whitney and Euler classes map to the
corresponding topological classes.
\end{prop}

Proof:
Both the algebraic and topological Stiefel-Whitney classes
(resp. Euler classes)
pull back to Chern classes of bundles on $Fl(V)$. Since $co$ maps
algebraic Chern classes to topological Chern classes
(\cite[Chapter 19]{Fulton}), the images of our algebraic
classes pull back on $Fl(V)$ to the pullbacks of the
corresponding topological classes. The proposition then follows
from  the naturality of $co$ and the
following lemma.
\begin{lemma} \label{thom.top}
If $f:Y \rightarrow X$ is a proper Chow envelope of smooth complex
varieties,
then the pullback $f^*:H^*(X;R) \rightarrow
H^*(Y;R)$ is injective for any coefficient group $R$.
\end{lemma}
Proof of Lemma \ref{thom.top}:
Consider the maps $cl:A_i(Y) \rightarrow H^{BM}_{2i}(Y)$ and
$cl:A_i(X) \rightarrow H^{BM}_{2i}(X)$
from Chow groups to Borel-Moore homology defined in
\cite[Chapter 19]{Fulton}. Let $s=[\tilde{X}] \in A_*(Y)$, where
$\tilde{X} \subset Y$ is a subvariety mapping birationally to $X$.
Then $f_*(s)=[X]$.
Since the class map is compatible with
proper pushforward,
$$f_*(cl(s))=cl([X]) = \mu_X,$$
where $\mu_X$ is the fundamental homology class of $X$.
If $x \in H^*(X)$, then by naturality of cap product,
$$f_*(f^*x \cap cl(s)) = x \cap f_*(cl(s)) = x \cap \mu_X.$$
On the other hand, $\_ \cap \mu_X$ is the Poincare
duality pairing (\cite[Chapter 19]{Fulton}). Thus, if
$f^*x = 0 \in H^*(Y;R)$, then the Poincare image of $x$ in
$H_*^{BM}(X;R)$ is zero. However, $X$ is smooth, so the pairing is
perfect, and thus $x=0$. Therefore $f^*$ is injective. \endproof

\section{Chow rings of quadric and flag bundles with
half integer coefficients} \label{s.halfchow} In this section we study
even rank bundles with quadratic form.  We will not assume that bundle
and its quadratic form are locally trivial in the Zariski topology.
Instead, we assume, throughout this section, that the structure group
reduces to $SO(2n)$.  We show (Theorem \ref{euler.haexist}) that
every such bundle $V \rightarrow X$ has characteristic classes $\pm
y_{n} \in A^{n}X$, such that $\frac{y_{n}}{2^{n-1}}$ is an Euler class
in $A^*(X; \Zh )$.  We use these Euler classes to compute the
Chow rings tensored with $\Zh$ of the associated quadric and
isotropic flag bundles.  There is no discussion of the odd rank case
because (over $\Z[\frac{1}{2}]$) the Chow groups of a quadric bundle
are generated by powers of the hyperplane section.

\paragraph{Construction of the isotropic flag bundle as a tower of
quadrics}

Let $V=V_n$ be a vector bundle of rank $2n$ with a non-degenerate
quadratic form.
Let $Q=Q_{n-1} \rightarrow X$  be the corresponding
quadric bundle. We have the following diagram:
$$\begin{array}{ccc}
Q & \stackrel{i}\rightarrow& \P(V)\\
& \searrow \scriptsize{\rho}  & \downarrow \scriptsize{\pi}\\
& & X
\end{array}$$

Let $V_{n}$ also denote the pullback bundle on $\P(V)$, $S_{n}$ the
tautological subbundle, and $S_{n}^{\perp}$ the orthogonal complement of
$S_{n}$ in $V_{n}$.  Although $S_{n}$ is not a subbundle of
$S_{n}^{\perp}$, the pullback $i^{*}S_{n}$ is a subbundle of
$i^{*}S_{n}^{\perp}$.
Set $V_{n-1} = i^*S_n^\perp/i^*S_n$; this is a rank $2n-2$ vector bundle
on $Q_{n-1}$ with nondegenerate
quadratic form.  Let $Q_{n-2} \rightarrow Q_{n-1} $  be the corresponding
quadric bundle.  In the same way, construct $V_{n-2}$, a rank
$2n-4$ vector bundle on $Q_{n-2}$.  Continuing, we get a tower
$$\begin{array}{cccccccccc}
Q_0 & \rightarrow & \P(V_1) & & & & & & &\\
&  & \downarrow  & & & & & &\\
& & Q_1 & \rightarrow & \P(V_2) & & & & &\\
& & & & \downarrow & & & & &\\
& & & & & \ldots & \rightarrow & \P(V_{n-1}) & &\\
& & & & & & & \downarrow & &\\
& & & & & & & Q_{n-1} & \rightarrow & \P(V_n)=\P(V)\\
& & & & & & & & & \downarrow\\
& & & & & & & & & X
\end{array}$$
(There is an analogous tower in the odd rank case, but it is not
essential to our discussion.)

Over $Q_{1}$ there is then a rank $2$ vector bundle $V_{1}$ with
nondegenerate quadratic form.  The corresponding quadric bundle $Q_{0}$
is a double cover of $Q_{1}$.

\begin{lemma}
\label{l:structuregroup}
With notation as above, if the structure group of $V_{n} \rightarrow X$ reduces
to $SO(2n)$, then the structure group of $V_{n-1} \rightarrow Q_{n-1}$ reduces
to $SO(2n-2)$. Moreover, if $V_n$ with its quadratic form is locally
trivial in Zariski topology, then $V_{n-1}$  (with its quadratic form)
is as well. \endproof

\end{lemma}
The proof is straightforward, and we omit it.

Since we are assuming that the structure group of $V_{n}
\rightarrow X$ reduces to $SO(2n)$, the above lemma
implies that the structure group
of $V_{1} \rightarrow Q_{1}$ reduces to $SO(2)$.

\begin{lemma} \label{l.rank2}
Let $V_1 \rightarrow X$ be a rank $2$ vector bundle with nondegenerate
quadratic form whose structure group reduces to $SO(2)$.  Then $V$
splits into a direct sum $V_1^\prime \oplus V_1^{\prime \prime}$ of
isotropic line bundles.
\end{lemma}
Proof: If $W$ is a vector space of rank 2 with a non-degenerate
quadratic form, then $W$ has exactly 2 isotropic lines. These
lines intersect in the origin, and each is invariant under
$SO(2)$. Thus, if the transition functions of the bundle $V_1$ lie in
$SO(2)$, we can construct two isotropic sub-linebundles of $V_1$
whose direct sum is $V$. \endproof

\medskip

To avoid cumbersome notation, we will temporarily write $V_i$, $S_i$,
etc. for the pullbacks of these bundles to $Q_1$.

The preceding lemma implies that the bundle $V_1 \rightarrow Q_1$
splits into a direct sum $V_1^\prime \oplus V_1^{\prime \prime}$.
This implies that each of the bundles $V_2, V_3, \ldots V_n$ has
two isotropic subbundles when pulled back to $Q_1$. They can be
defined as follows. Assume by induction that $V_{i-1}^\prime$ and
$V_{i-1}^{\prime\prime}$,
maximal isotropic subbundles of $V_i$ are defined. By definition
we have a surjective map
$$S_i^\perp \rightarrow S_i^\perp/S_i = V_{i-1}.$$
Define $V_i^\prime$ (resp. $V_i^{\prime \prime}$) to be the preimage of
$V^\prime_{i-1}$ (resp $V_{i-1}^{\prime \prime}$). Then
$V_i^\prime$ and $V_i^{\prime \prime}$ are maximal isotropic subbundles
of $V_i$. Note that by construction,
\begin{equation} \label{eqn} c_i(V_i^\prime)  = h_i c_{i-1}(V^\prime_{i-1})
\end{equation} where
$h_i = c_1(S_i^*)$ is the hyperplane class on  $\P(V_i)$.

The space $Q_1$ can be identified with the bundle
$Fl(V)$ of isotropic flags in $V$ of length $n-1$. On
$Fl(V)$ the bundle $V$ has a tautological flag of
subbundles, which we will denote $E_1 \subset E_2 \subset \ldots
\subset E_{n-1}$. In this notation, the bundles $V_i$ (pulled back to
$Fl_{n-1}$) are the quotients $E_{n-i}^\perp/E_{n-i}$ (taking $E_0 =
0$).  The bundle $V$ has two maximal isotropic subbundles
$V_n^\prime$ and $V_n^{\prime \prime}$. Using the fact that the
bundles $V_1^\prime$ and $V_1^{\prime \prime}$ are dual, direct
calculation, via Equation \ref{eqn},
shows that $$c_n(V_n^\prime) = -c_n(V_n^{\prime \prime}).$$

\begin{prop} \label{pushpull}Let $f:Fl(V_n) \rightarrow X$ be the projection.
Then
there is a (canonical) class $s \in A^*(Fl(V_n))$ such that for any $Y
\rightarrow X$ and any
$x \in A^*(X)$, $f_*(s \cdot f^*x \cap [Fl(V_n) \times_{X}Y ] ) =
2^{n-1}x \cap [Y]$.  Hence $f^*$ is injective with $\Zh$ coefficients.
\end{prop}
Proof: We identify $Fl(V_n) = Q_1$. If  $\rho_i:Q_{i-1} \rightarrow Q_{i}$ is
a quadric in our tower, $x_{i} \in A^*(Q_i)$, and $Y \rightarrow
Q_{i}$,
then $\rho_{i*}(h_i^{2i-2} \cdot  \rho_i^* x \cap [Q_{i-1} \times_{X}Y] )
= 2x \cap [Y]$.
Thus $s = h_2^2 h_3^4 \ldots h_{n}^{2n-2}$ (pulled back to $Fl(V_n))$ is the
desired class. \endproof

The following is an immediate consequence of Proposition \ref{pushpull}.
\begin{thm} \label{euler.haexist}
Let $V$ be a vector bundle of rank $2n$ with a non-degenerate quadratic
form (not necessarily locally trivial in the Zariski topology)
whose structure
group reduces to $SO(2n)$.
Then there are characteristic classes $\pm y_n \in A^n(X)$ such that
$\pm\frac{y_n}{2^{n-1}}$ are Euler classes in $A^n(X; \Zh )$.
\end{thm}
{\bf Remark:} Theorem \ref{euler.haexist} strengthens a result
of Vistoli \cite{Vistoli} who showed (implicitly) the existence
of an Euler class in $A^*X \otimes \Q$ for arbitrary principal
$SO(2n)$ bundles.
\medskip

Proof: The existence of the class $s$
from Proposition \ref{pushpull} implies that if
$c$ is in the kernel of the map $A^{*} Fl(V)  \stackrel{p_1^* - p_2^*}
\rightarrow Fl(V) \times_X Fl(V)$
then the class $d$ in $A^*X$ defined by the formula
$$d \cap [Y] = f_*(s \cdot c \cap [Fl(V_n) \times_{X}Y ])$$
satisfies $f^*d = 2^{n-1} c$.  This is an
analogue of \cite[Theorem 2.3]{Kimura} used above,
and is proved in the same way.  Setting  $c=c_n(E)$
where $E$ is one of the two tautological maximal
isotropic subbundles of $f^*V$, produces a class $y_n \in A^{*} X$ such
that $f^*y_n = 2^{n-1}c_n(E)$.  The classes $\pm y_n$ are natural
with respect to pullbacks, because the class $s$ of Proposition
\ref{pushpull} is natural.
\endproof

\paragraph{Computation of Chow groups}

We now describe $A^*(Q; \Zh)$
as an algebra over $A^*(X ;\Zh)$.
Iterating over the quadric tower used above, we are
also able to compute $A^*(Fl(V);\Zh)$ as an algebra over
$A^*(X; \Zh )$.

Before we state the
theorems, we need some further properties of Euler classes.

\begin{prop} \label{euler.recur}
Let $x_{n-1}$ be an Euler class for the associated
bundle $V_{n-1} \rightarrow Q_{n-1}$, and let $h$ be the pullback
of the hyperplane section on $\P(V)=\P(V_n)$. Then $hx_{n-1} = \rho^*x_n$,
where $x_n \in A^n(X;\Zh)$ is an Euler class for $V$. If
the bundle is locally trivial, then the identity holds integrally.
\end{prop}
Proof: Let $f:Fl(V) \stackrel{g} \rightarrow Q \stackrel{\rho} \rightarrow X$
be the projection. Then
$f^*x_n = c_n(V^\prime_n)$ where $V_n^\prime$
is a maximal isotropic subbundle of $f^*V$.
On the other hand, $g^*x_{n-1} = c_{n-1}(V_{n-1}^\prime)$ where
$V_{n-1}^\prime$ is a maximal isotropic subbundle of $f^*V_{n-1}$.
By Equation \ref{eqn}, $c_n(V_{n}^\prime) =hc_{n-1}
(V_{n-1}^\prime)$, so the relation holds on $Fl(V)$. Since
$f^*$ is injective up to 2-torsion (or in the locally trivial
case -- integrally injective) the relation follows.
\endproof

\medskip

The following theorem is the key connection between Euler classes
and Chow groups of quadric bundles. We use the same notation as
Proposition \ref{euler.recur} above.
\begin{thm} \label{euler.fiber}
On a fiber of $Q_{n-1} \rightarrow X$ the Euler class $x_{n-1}$ restricts to
$\pm (e_{n-1} -f_{n-1})$,
where $e_{n-1}$ and $f_{n-1}$ are the two ruling classes on the fiber.
\end{thm}
Proof: Since we are working on a fiber, we may assume $X$ is a point,
and $V=V_n$ is a vector space.
Consider the tower above. The space $Q_0$ is the flag variety of
length $n$ isotropic flags in $V$.
It is a disconnected double cover of
$Q_1 = Fl(V)$. Any element $g \in O(2n)$ induces an
automorphism of this tower. At each step denote
the induced maps $Q_i \rightarrow Q_i$ and $\P(V_i) \rightarrow \P(V_i)$ by
$g_i$. By construction, the $g_i$'s are compatible for different values of
$i$. If $g$ is not in the identity component of $O(2n)$, then
$g_0:Q_0 \rightarrow Q_0$ acts by exchanging the two sheets of $Q_0$.

\begin{lemma} \label{euler.auto}
Let $g$ be an element of $O(2n)$ not in the identity component. Then
$g_1^*(c_i(V_i^\prime) = -c_i(V_i^\prime)$.
\end{lemma}
Proof of Lemma \ref{euler.auto}:
It is clear that the hyperplane class in any of the projective
bundles in the tower is invariant under pullback by the automorphism
induced by $g$. Since $c_i(V_i^\prime) = h_ic_{i-1}(V_{i-1}^\prime)$,
it suffices to prove the proposition when $i=1$. Now $Q_0$ is
the union of two components, $\P(V_1^\prime)$ and $\P(V_1^{\prime\prime})$,
and
$i_{0_*}([\P(V_1^{\prime\prime})] - [\P(V_1^{\prime})]) = \pi_1^*(2c_1
(V_1^\prime)).$

Since $g$ exchanges the two sheets of $Q_0$,
$$g_0^*([\P(V_1^\prime)] -[\P(V_1^{\prime\prime})]) = -([\P(V_1^\prime)] -
[\P(V_1^{\prime\prime})]).$$ Thus,
$$g_0^*(\pi_1^*(2c_1(V_1^\prime)))= \pi_1^*(g_1^*(2c_1(V_1^\prime)))
=-\pi_1^*(2c_1(V_1^\prime)).$$ Since $\pi_1^*$ is injective,  and
we are working with coefficients in $\Z[\frac{1}{2}]$, the lemma
follows. \endproof

\medskip

Since
$e_{n-1}$ and $f_{n-1}$ generate $A_{n-1}(Q_{n-1})$ (Lemma
\ref{quad}), we can write $x_{n-1} = ae_{n-1} + bf_{n-1}$, where $a$
and $b$ are in $\Z[\frac{1}{2}]$. Thus, $i_*(x_{n-1}) =(a+b)h_n$.
Let $g$ be an element of $O(2n)$ not in the identity component.  Then
in the notation of the preceding proposition, $g_{n-1}^*x_{n-1}=
-x_{n-1}$. Thus $a+b =0$.
Hence we can write $x_{n-1} = \lambda(e_{n-1} - f_{n-1})$
for some $\lambda\in \Z[\frac{1}{2}]$.
Since $x_{n-1}$ is an Euler class, $x_{n-1}^2=(-1)^nc_{2n-2}(V_{n-1})$
(Proposition \ref{p.basic}).  By direct calculation,
$$c(V_{n-1}) = \frac{c(V_n)}{1-h_{n}^2}.$$
But $c(V_n) = 1$ since $V_n$ is trivial.  Thus,
$$ x_{n-1}^2= (-1)^{n} c_{2n-2}(V_{n-1}) = (-1)^{n} h_n^{2n-2}=(-1)^{n} 2.$$
On the other hand, $x_{n-1}^2=\lambda^2(e_n-f_n)^2$, so by Lemma \ref{rels},
$x_{n-1}^2 = \lambda^2 \cdot (-1)^n\cdot 2$. Therefore, $\lambda = \pm 1$
as desired. \endproof

\begin{thm} \label{q.algebra} Let $X$ be a scheme and
let $V\rightarrow X$ be a vector bundle of rank $2n$
whose structure group reduces to $SO(2n)$ with respect to a non-degenerate
quadratic form on $V$. Let $Q \subset \P(V)$ be the associated quadric
bundle.
Then
$A^*(Q;\Z[\frac{1}{2}])= A^*(X;\Zh)[h,x_{n-1}]/I $, where
$I$ is the ideal generated by the relations $$hx_{n-1} = \rho^*(x_n)$$
$$x_{n-1}^2 = (-1)^{n-1} c_{2n-2}(V_{n-1}) = (-1)^{n-1}(h^{2n-2} +
h^{2n-4}c_{2}(V) + \ldots + h^{2}c_{2n-4}(V) + c_{2n-2}(V)).$$
\end{thm}

{\bf Remark:} The relation $$h^{2n} + h^{2n-2}c_{2}(V)+ \ldots
+ c_{2n}(V) = 0 ,$$ inherited from $A^*(\P(V))/A^*(X)$, can
be easily derived from the two relations above.
If $V \rightarrow X$
has rank $2n+1$ then
$$A^*Q=A^*X[h]/h^{2n+1} + h^{2n-1}c_2(V) + \ldots + c_{2n}(V).$$

{\bf Example:} Theorem \ref{q.algebra} can be used to show that not
all rank vector bundles of rank $N$ whose structure group reduces to
$SO(N)$ have maximal isotropic subbundles.  For an example when $N$ is
even, let $Q_2 \subset \P^5$ be a smooth quadric of rank 3. On $Q_2$
we can define the bundle $V_2=S^\perp/S$, where $S={\cal O}(1)$, and
$S^\perp$ is the orthogonal bundle with respect to the quadratic form
on $Q$. Direct calculation shows that $V_2$ is locally trivial in the
Zariski topology, so the bundle has rank 4 and structure group reduces
to $SO(4)$. We claim, however, that $V_2$ has no isotropic subbundle
of rank 2. Suppose to the contrary, that $E \subset V_2$ was such a
bundle. Then $c_2(E) = x_2$, while $c_1(E) = ah$ for some constant
$a$, since $h$ generates $A^1(Q)$.  On the other hand,
$$c(E)c(E^*)=c(V_2)=1+h^2+h^4.$$ Thus, $$(1
+ah+x_2)(1-ah+x_2)=1+h^2+h^4.$$ In particular, $a^2h^2+2x_2=h^2$,
which is impossible since $h^2$ and $x_2$ are independent in $A^2(Q)$.
Therefore, no such bundle $E$ exists.  For an example of an odd rank
bundle, let $V$ be the rank 3 bundle $S^\perp/S$ on the quadric
hypersurface in $\P^4$. A
similar calculation (using only the hyperplane class) shows that $V$ has
no isotropic line subbundles.

\medskip

Now let $f:Fl(V) \rightarrow X$ be the
flag bundle  over $X$. Let $E_1 \subset E_2 \subset \ldots \subset
E_{n-1} \subset E_{n} \subset f^*V$ be one of the two tautological isotropic
flags. Via the quadratic form, we can extend it
to a complete flag
$E_1 \subset \ldots \subset E_n \ldots \subset E_{2n} =f^*V$.
Set $h_i=c_1(E_{n-i+1}/E_{n-i})$ (where $E_0=0$).

\begin{thm} \label{fl.algebra} Let $V$ be as in the statement of Theorem
\ref{q.algebra}. Then
$$A^*(Fl(V); \Zh ) = A^*(X; \Zh)[h_1, \ldots, h_n]/I $$
where $I$ is the ideal generated the relations
$$\prod_{i=1}^{n} (1-h_i^2) = f^*c(V)$$
$$h_1\;h_2 \ldots h_n=f^*x_n$$
\end{thm}
{\bf Remark:} The choice of notation $h_i$ in the above theorem
is consistent with its use previously, since for $i>1$,
$c_1(E_{n-i+1}/E_{n-i})$
is the pullback of the hyperplane class of the projective bundle
$\P(V_i)$.\\

Proofs of Theorems \ref{q.algebra} and \ref{fl.algebra}:
We give the proof assuming $X$ is smooth to avoid
using the language of bivariant intersection theory
(\cite[Chapter 17]{Fulton}).
A reader familiar with this theory can easily extend the proof
to singular schemes.

Theorem \ref{fl.algebra} follows from Theorem \ref{q.algebra} by
iterating over the tower of quadrics constructed previously.
Thus it suffices to prove Theorem \ref{q.algebra}.
First note that $A^*X[h,x_n]/I$ has a basis of monomials
$h^ix_n^\alpha$,
where $0 \leq i \leq n-1$, and $\alpha \in \{0,1\}$. Since
$x_{n-1}$ restricts on a fiber to the difference of two rulings,
these classes restrict to a basis for the Chow group tensored with $\Zh$
of a fiber (Lemma \ref{quad}). At this point we would like to
apply Noetherian induction and conclude that $h^ix_n^\alpha$ is
a basis for the Chow groups of $Q$ over $A^*(X)$. Unfortunately, the bundle
$Q \stackrel{\rho}  \rightarrow X$ need not be locally trivial in
the Zariski topology. However,
$Q \times_X Fl(V)$ is Zariski
locally trivial over $Fl(V)$ since the pullback of $V$ has
a maximal isotropic subbundle.
Thus, by Noetherian induction,
our classes pull back to
a basis for $A^*(Q\times_X Fl(V); \Zh )$ over $A^*(Fl(V))$.
On the other hand, the pullback
$f^*:A^*(X; \Zh ) \rightarrow A^*(Fl(V); \Zh )$ has a
section (see Proposition \ref{pushpull}). Pushing forward via the
section, we conclude that $\{h^ix_n^\alpha\}$ is a basis for
 $A^*(Q; \Zh )$.

To complete the proof, we must check the relations.
The relation $hx_{n-1}= \rho^*(x_n)$ is Proposition \ref{euler.recur}.
Since $x_{n-1}$ is an Euler class for the bundle $V_{n-1}$,
$x_{n-1}^2 = (-1)^{n-1}c_{2n-2}(V_{n-1})$ (Proposition \ref{p.basic}).
Since $c_{2n-2}(V_{n-1}) = \left\{ \frac{c(V)}{1-h^2}
\right\}_{2n-2}$, the second relation follows.
\endproof

\medskip

{\bf Remark:} In general $A^*(Fl(V))$ need not be
free over $A^*(X)$. If it were, then one could define
an integral Euler class for bundles whose structure
group reduces to $SO(2n)$, but were not locally trivial
in the Zariski topology. However, we noted before that B. Totaro
(\cite{T1}) has shown that this can not be true.

{\bf Remark:}
If the stucture group $V \rightarrow X$ does not reduce to
$SO(2n)$, then the flag variety $f:Fl_n(V) \rightarrow X$
of length $n$ isotropic
flags is connected. If
$E \subset f^*V$ is a tautological rank $n$ isotropic subbundle, then
arguments similar to those used in Lemma \ref{euler.auto}
shows that $c_n(E) = 0 \in A^*(Fl_n(V); \Zh )$. Thus
the Euler class is also 2-torsion, a fact consistent
with topology, since non-orientable bundles only have
Thom classes with $\Z/2\Z$ coefficients.
This calculation also shows that $c_{2n}(V)$ is 2-torsion
as well.
\medskip
\section{Chow groups of quadric and isotropic flag bundles with integer
coefficients} \label{s.zchow}
Let $V \rightarrow X$ be a vector bundle with a quadratic form, and let
$F \subset V$ be a fixed maximal isotropic subbundle. The purpose
of this section is to compute the
Chow rings of the associated quadric and flag bundles. The presentation
will depend on the particular subbundle $F$.

\paragraph{Chow rings of quadric bundles}
\begin{thm} \label{zch.quad}
Let $V$ be a vector bundle of rank $N$ with
with a nondegenerate
quadratic form on a scheme $X$.
Let $Q \stackrel{\rho} \rightarrow X$
be the associated quadric bundle and $h$ the hyperplane
class on $Q$.
Assume that $F \subset V$ is a maximal isotropic subbundle
and set $\gamma = [\P(F)] \subset Q$. Then
$$A^*(Q) = A^*X[h,\gamma]/I$$ where
$I$ is the ideal generated by the relations
$$2h\gamma= h^{n} -c_1(F)h^{n-1} + \ldots + (-1)^nc_n(F)$$
$$\gamma^2= (-1)^{n-1}(c_{n-1}(F) + c_{n-3}(F)h^2 + \ldots )\gamma$$
if $N=2n$, and the relations
$$2h\gamma = h^{n+1} + c_1(V/F)h^{n} + \ldots c_{n+1}(V/F)$$
$$\gamma^2 = (c_n(V/F) +c_{n-2}(V/F)h^2 + \ldots)\gamma$$
if $N=2n+1$.
\end{thm}
Proof: We give the proof only in the even case, as the odd case is analogous.
As in the proof of Theorem \ref{euler.haexist} we assume $X$ is smooth.
A basis for $A^*X[h,\gamma]/I$  as an $A^*X$ module
is the monomials
$1,h,h^2, \ldots h^{n-1},\gamma,h\gamma, \ldots , h^{n-1}\gamma$. These
classes restrict to a basis of the Chow groups
of the fiber, so they form a basis for $A^*(Fl(V))$ as a $A^*(X)$ module
(All bundles are locally trivial in Zariski topology, so we can apply
Noetherian induction).
To prove the theorem we must check the relations.

Let $i:Q \hookrightarrow \P(V)$ be the inclusion. From the proof of
Theorem \ref{t.chern}, we know that
$i_*\gamma = H^n -c_1(F)H^{n-1} + \ldots + (-1)^n c_n(F)$, where
$H$ is the hyperplane section on $\P(V)$. Hence the right hand side
of the first relation is $i^*i_*\gamma$. On the other hand,
the normal bundle to $Q$ in $\P(V)$ is $2h$. Thus, by the self
intersection  formula, $i^*i_*\gamma= 2h\gamma$.
This proves the first relation.

The second relation follows from the self intersection formula
$\gamma^2=j_*(c_{n-1}(N_{\P(F)}Q))\gamma$ where $j_*:\P(F) \hookrightarrow
Q$ is the inclusion. The normal bundle is
$$N_{\P(F)}Q=N_{\P(F)}\P(V)/N_{Q}\P(V) = \frac{V/F \otimes S^*}{(S^*)^{\otimes
2}}.$$
Tensoring the righthand side with $S \otimes S^*$, and identifying $V/F$
with $F^*$
we can identify the normal bundle as $\frac{F^*}{S^*}\otimes S^*$. The formula
then follows by direct calculation.
\endproof

\paragraph{Splitting bundles with quadratic form}
Let $V= L_1 \oplus L_2 \oplus \ldots \oplus L_N$ be a direct sum of
line bundles of the form $L_1 \oplus L_2 \oplus \ldots \oplus L_n
\oplus L_n^* \oplus \ldots \oplus L_1^*$ if $N=2n$, and
$L_1 \oplus \ldots \oplus L_n \oplus M \oplus L_n^* \oplus \ldots
\oplus L_1$, if $N=2n+1$, where $M$ is a self-dual line bundle.
Define a quadratic form $q_{std}$ on $V$ by the rule that
$q_{std}$ restricted to $L_i \oplus L_{N-i+1}$ is the canonical
pairing, and $L_i$ is orthogonal to $L_j$ otherwise. We will
call $q_{std}$ the standard quadratic form on $V$.
Suppose that $q$ is a quadratic form on $V$ such that
$q$ restricted to $L_i \oplus L_{N-i+1}$ is the canonical pairing
and such that (with respect to $q$), $L_i$ is orthogonal to
$L_1, \ldots , L_{N-i}$.

\begin{lemma}
Let $V$, $q$ and $q_{std}$ be as above. Then, (assuming
that the base field has characteristic not equal to 2),
there exists a
vector bundle automorphism $g$ of $V$ such that $g$ is the identity
on $L_1 \oplus \ldots L_n$ and such that $g^*q_{std} = q$.
\endproof\end{lemma}

The proof of this lemma is an exercise in ``bundle-ized'' linear
algebra which we omit.

{\bf Remark:} Fulton \cite{F1} shows
that there is a deformation of $(V,q)$ to $(V,q_{std})$
that induces an isomorphism of the Chow rings of the associated
isotropic flag bundles. As a result, Theorems \ref{ch.dn}
and \ref{ch.bn} are still valid in characteristic 2.

The following quadratic
splitting principle (in characteristic not equal to 2)
is an application of the lemma above.
\begin{prop} \label{q.split}
Let $V \rightarrow X$ be a vector bundle of rank $N=2n$ or $N=2n+1$
with a nondegenerate quadratic
form. Let $E_1 \subset \ldots E_n$ be a
maximal isotropic flag of subbundles of $V$. Set $L_i=E_i/E_{i-1}$
and if $N$ is odd, set $M=E_n^\perp/E_n$. Then there exists
$f:Y \rightarrow X$ with $f^*:A^*X \rightarrow A^*Y$ injective, and
an isomorphism of $f^*V$ with $L_1 \oplus \ldots L_n \oplus L_n^*
\oplus \ldots \oplus L_1^*$ or
$L_1 \oplus \ldots L_n \oplus M \oplus L_n^*
\ldots \oplus L_1^*$, as vector bundles with quadratic form (
where $f^*V$ is given the quadratic form inherited from $V$, and
the direct sum of line bundles is given the standard form described
above). Furthermore, the subbundle $E_i$ corresponds to the sum
$L_1 \oplus \ldots \oplus L_i$.
\end{prop}
Proof:
Extend $E_1 \subset \ldots \subset E_n$ to a complete flag
$E_1 \subset \ldots \subset E_N = V$
by setting $E_{n+i}=E_{n-i}^\perp$. Then
$L_{n+i} =E_{n+i}/E_{n+i-1}
\simeq L_{n-i+1}^*$. Then the complete splitting principle as
stated in \cite{F1} (cf. \cite[Theorem 8.3]{Gillet}) gives
$f:Y \rightarrow X$ with $f^*$ injective and and isomorphism
of
$f^*V$ with
$L_1 \oplus \ldots L_n \oplus L_n^*
\oplus L_1^*$ or $L_1 \oplus \ldots L_n \oplus M \oplus L_n^*
\ldots \oplus L_1^*$, as vector bundles. Let $q$ be the quadratic form
on $L_1 \oplus \ldots \oplus L_N$ pulled back by isomorphism from
the quadratic form on $f^*V$. By construction, $q$ restricts
to the dual pairing on $L_i \oplus L_{N-i+1}$. Moreover, since we began
with a maximal isotropic flag, each $L_i$ is orthogonal to
$L_1, \ldots L_{N-i}$. By the above lemma, there exists a vector
bundle automorphism of
$L_1 \oplus\ldots \oplus L_N$ taking $q$ to $q_{std}$.
This proves the proposition.
\endproof
\paragraph{Presentation of the Chow rings of flag bundles}
Let $V$ be a vector bundle
with a nondegenerate quadratic form on a scheme $X$.
We assume the existence of a maximal isotropic flag
$F_{\cdot}= F_1 \subset F_2 \ldots \subset F_n \subset V$.
To compute the Chow
ring we need only assume the existence of a maximal isotropic subbundle
$F$. However, assuming the existence of the full flag $F_\cdot$ gives
a presentation consistent with the Schubert variety ideas of \cite{F1}.

Let $\pi:Fl(V) \rightarrow X$ be the bundle of
maximal isotropic flags in $V$. (Recall that if the rank
$V$ is $2n$, then the flag variety has two isomorphic connected components.
In this case the notation $Fl(V)$ will refer to one of these
components, which may be identified with $Fl_{n-1}(V)$.)
Set $y_i = -c_1(\pi^*F_i/\pi^*F_{i-1})$. The following two
theorems (depending on whether the rank of $V$ is even or odd)
describe $A^*Fl(V)$ as an algebra over $A^*X$.

\begin{thm} \label{ch.dn}
({\em $D_n$ case}): If $V$ has rank $2n$ then
$A^*Fl(V) =A^*X[x_1,x_2, \ldots , x_n,c_1, \ldots c_{n-1}]/I$
where $I$ is the ideal generated by the following relations (for
$1 \leq i \leq n$):
$$e_i(-x_1^2,\ldots ,-x_n^2) = \pi^*c_{2i}(V)$$
$$(-1)^pc_p^2 + (-1)^{p-1}2c_{p-1}c_{p+1} + \ldots - 2c_{2p-1}c_{1}
+2c_{p} = c_{2p} - c_{2p-1}e_1(x_1,\ldots x_n) + \ldots + (-1)^pc_p
e_p(x_1, \ldots , x_n)$$
$$2c_i = e_i(x_1, \ldots , x_n) + e_i(y_1, \ldots , y_n)$$
($c_k=0$ for $k \geq n$).
\end{thm}

Recall that if $a_1, \ldots a_n$ is a set of variables then
$e_i(a_1, \ldots , a_n)$ denotes the $i$-th elementary symmetric
polynomial in those variables.

\begin{thm} \label{ch.bn}
({\em $B_n$ case}): If $V$ has rank $2n+1$, then
$A^*Fl(V) = A^*X[x_1, \ldots , x_n,c_1, \ldots , c_n]/I$
where $I$ is the ideal generated by the following relations (for
$1 \leq i \leq n$)
$$e_i(-x_1^2, \ldots , -x_n^2) = \pi^*c_{2i}(V)$$
$$le_i(-x_1^2, \ldots , -x_n^2)  = \pi^*c_{2i+1}(V)$$
$$2c_i  =  e_i(x_1+l, \ldots , x_n +l) + e_i(y_1+l,y_2+l, \ldots , y_n +l)$$
\begin{eqnarray*}
(-1)^pc_p^2 & + & (-1)^{p-1}2c_{p-1}c_{p+1} + \ldots - 2c_{2p-1}c_{1}
+2c_{p}\\
&=& c_{2p} - c_{2p-1}e_1(x_1+l,\ldots x_n+l) + \ldots + (-1)^pc_p
e_p(x_1+l, \ldots , x_n+l)
\end{eqnarray*}

where $l = \pi^*c_1(V)$
\end{thm}

Proof of Theorem \ref{ch.dn}.
If $V$ has rank $2n$, let $E_1 \subset E_2 \ldots \subset E_n$
be a tautological flag on $\pi^*V$ such that $E_n$ and $\pi^*F_n$
are in opposite ruling families. Set $x_i = c_i(E_i/E_{i-1})$.
Since $\pi^*V/E_n = E^*_n$,  $e_i(-x_1^2, \ldots -x_n^2)= c_{2i}(\pi^*V)$
proving the first relation.

Let $\rho:Q \rightarrow Fl(V)$ be the quadric bundle associated
to $\pi^*V$. Define classes $c_i \in A^*(Fl(V))$ by the relation
$$[\P(E_n)] + [\P(\pi^*F_n)] = h^{n-1} - h^{n-2}\rho^*c_1 +
h^{n-3}\rho^*c_2 - \ldots + (-1)^{n-1}\rho^*c_{n-1}$$
Since $\rho_* (h^n[\P(E_n)]\rho^*c_i)
= c_i$, the $c_i$'s are uniquely determined. Furthermore, by the
proof of Theorem \ref{t.chern}, $2c_i = c_i(E_n) + c_i(F_n)
=e_i(x_1, \ldots, x_n) + e_i(y_1, \ldots, y_n)$.

Note that $A^*X[x_1, \ldots x_n,c_1, \ldots c_{n-1}]/I$ has a basis
of monomials of the form
$$x_1^{a_1}x_2^{a_2}\cdots x_n^{a_n}
c_1^{\alpha_1}c_2^{\alpha_2} \cdots c_{n-1}^{\alpha_{n-1}}$$
where $0 \leq a_i \leq n-i$ and $\alpha_i \in \{0,1\}$.
On the other hand, Fulton \cite{F1} (generalizing a result stated
in \cite{Marlin}) has shown that these
monomials form a basis for $A^*(Fl(V))$ as an $A^*X$ module.
Therefore, to prove the theorem it suffices to check that
the second relation holds.

By Proposition \ref{q.split} there is a map
$f:Y \rightarrow Fl(V)$ such that $f^*$ is injective and
$\pi^*V$ pulls back to $L_1 \oplus L_2 \ldots L_n \oplus L_n^*
\oplus \ldots \oplus L_1$ with the standard quadratic form.
To check a relation, it suffices to check it for the totally
split bundle
on $Y$.

Furthermore, there is always a map $g:Y^\prime \rightarrow Y$
with $g^*$ injective such that $Y^\prime$ is quasi-projective
(this is Chow's Lemma combined with Nagata's
embedding theorem, see \cite[Section 18.3]{Fulton}).
We may therefore assume that $Y$ is quasi-projective.
In particular
there are  maps $j_i:Y \rightarrow \P^{k_i}$
such that $L_i = j_i^*M_i$ for some line bundle $M_i$ on
$\P^{k_i}$. Setting $Z=\P^{k_1} \times \P^{k_2} \times \ldots
\times \P^{k_n}$, there is a map $j:Y \rightarrow Z$
such $L_1 \oplus L_2 \oplus \ldots L_n \oplus L_n^* \oplus \ldots \oplus
L_n^*$ is the pullback of direct sum of line bundles and their
duals on $Z$. To prove the relation, it suffices to check it in
$A^*Z$. However, since $A^*Z$ is torsion free, we need only check that
the relation holds up to multiplication by a constant.
Since $2c_p = e_p(x_1, \ldots , x_n) + e_p(y_1, \ldots y_n)$,
multiplying both sides of the potential relation by 4 we reduce to showing
that the degree $2p$ term in
$$(c(E_n) + c(\pi^*F_n))(c(E_n^*) + c(\pi^*F^*_n))$$ equals the
degree $2p$ term in
$$c(E_n)(c(E^*_n) + c(\pi^*F_n)) + c(E_n^*)(c(E_n) + c(\pi^*F_n)).$$
This equality is immediate since
$$c(E_n)c(E_n^*) = c(\pi^*F_n)c(\pi^*F_n^*) = c(\pi^*V).$$ This proves
Theorem \ref{ch.dn}.
\endproof

Proof of Theorem \ref{ch.bn}.
If $V$ has rank $2n+1$, let $E_1 \subset E_2 \subset \ldots  \subset E_n$
be the tautological isotropic flag on $Fl(V)$. Again set
$x_i = c_1(E_i/E_{i-1})$. Let $L=E_{n+1}/E_{n}$. Then
$l = c_1(\pi^*V)$, and $c(V)=c(E)c(E^*)(1+l)$. This shows that the
first two relations hold.

Again let $\rho:Q \rightarrow Fl(V)$ be the quadric bundle
associated to $\pi^*(V) \otimes L$. Define the $c_i$ by the equation
$$[\P(E_n \otimes L) + \P(\pi^*F_n \otimes L)] =
h^n -h^{n-1}\rho^*c_1 +h^{n-1}\rho^*c_2-
\ldots + (-1)^n \rho^*c_n.$$ The proof of Theorem \ref{t.chern}
also shows that $2c_i = c_i(E_n \otimes L) + c_i(\pi^*F_n \otimes L)$,
thereby proving the last relation.

Fulton \cite{F1} has also shown that the appropriate monomials
generate $A^*(Fl(V))$ as an $A^*(X)$ module, so that proving the
theorem reduces to checking the quadratic relation on the $c_i$'s.

The proof of the quadratic relation in the $B_n$
case is more or less the same
as that given above. There is one twist. The totally split pullback of
$\pi^*V$ is $L_1 \oplus L_2 \oplus \ldots L_n \oplus M \oplus L_n^*
\oplus \ldots \oplus L_1^*$, where $M$ is self dual.
However, $\pi^*(V)\otimes M$ totally splits
as $M_1 \oplus M_2 \ldots \oplus M_n  \oplus {\cal O}
\oplus M_n^* \oplus \ldots \oplus M_1^*$ (where $M_i = L_i \otimes M$).
We can then
reduce to checking the relation on a product of projective spaces.
The relation
on the $c_p$'s then follows from the fact that
$$c(E_n \otimes M)c(E^*_n \otimes M) = c(\pi^*F_n \otimes M)
c(\pi^*F_n^* \otimes M).$$
This proves Theorem \ref{ch.bn}.
\endproof



\begin{thebibliography}{99}
\bibitem[EG]{E-G-T} D. Edidin, W. Graham, {\it
Characteristic classes in the Chow ring}
preprint.
\bibitem[F1]{F1} W. Fulton, {\it Schubert varieties in flag bundles},
preprint 1994.
\bibitem[F2]{F2} W. Fulton, {\it Determinantal formulas for
orthogonal and symplectic degeneracy loci}, preprint 1994.
\bibitem[Fulton]{Fulton} W. Fulton, {\it Intersection Theory}, Ergebnisse, 3.
Folge, Band 2,
Springer Verlag, (1984).
\bibitem[G]{Gillet} H. Gillet, {\it Riemann-Roch theorems for higher
$K$-theory}, Adv. Math. {\bf 40} (1981), 203-289.
\bibitem[K]{Kimura} S. Kimura, {\it Fractional intersection theory
and bivariant theory}, Comm. Alg. {\bf 20} (1992) 285-302.
\bibitem[M]{Marlin} R. Marlin, {\it Anneaux de Chow des groupes
algebraiques $SU(n),Sp(n),SO(n)$, $Spin(n),G_2,F_4$; torsion},
C.R. Acad Sc. Paris, Series A, {\bf 279} (1974), 119-122.
\bibitem[Milnor-Stasheff]{Milnor-Stasheff}
J. Milnor, J. Stasheff, {\it Characteristic Classes},
Princeton University Press (1974).
\bibitem[T1]{T1} B. Totaro, preprint 1994.
\bibitem[T2]{T2} B. Totaro, preprint 1994.
\bibitem[Sem-Chev]{Sem-Chev} Anneau de Chow et applications, Seminaire
Chevalley, Secr\'etariat math\'ematique, Paris (1958).
\bibitem[S]{Swan} R. Swan {\it }, {\it $K$-Theory of quadric
hypersurfaces},  Annals of  Math {\bf 122 } (1985), 113-153.
\bibitem[V]{Vistoli} A. Vistoli,
{\it Characteristic classes of principal bundles in algebraic intersection
theory}, Duke Math J. {\bf 58}
(1989), 299-315.
\end{thebibliography}
\end{document}